\documentclass[11pt]{amsart}

\usepackage[usenames]{color}
\usepackage{amssymb, amscd, amsmath, float, graphicx, longtable,
  tabularx, stmaryrd, fullpage, mathtools,todonotes, empheq }
\usepackage{amsfonts,amsxtra}
\usepackage{enumerate,enumitem,epsfig,verbatim}
\usepackage{fouriernc}\DeclareMathAlphabet{\mathcal}{OMS}{cmsy}{m}{n}

\usepackage[all,2cell]{xy}
\SelectTips{cm}{10} 
\usepackage{url}
\usepackage{datetime}

\usepackage[]{hyperref}

\theoremstyle{plain} 
\newtheorem{thm}{Theorem}[section]

\newtheorem{cor}[thm]{Corollary}
\newtheorem{lem}[thm]{Lemma}

\theoremstyle{definition}
\newtheorem{dfn}[thm]{Definition}

\newtheorem{conj}[thm]{Conjecture}
\newtheorem{ques}[thm]{Question}

\newtheorem{conv-not}[thm]{Conventions and Notation}

\newcommand{\sk}{\mathsf{k}}
\newcommand{\m}{\mathfrak{m}}
\newcommand{\n}{\mathfrak{n}}
\newcommand{\p}{\mathfrak{p}}

\newcommand{\ff}{\mathfrak{f}}

\newcommand{\ZZ}{\mathbb{Z}}

\renewcommand{\to}{\longrightarrow}

\DeclareMathOperator{\rank}{rank}

\DeclareMathOperator{\Ext}{Ext}

\DeclareMathOperator{\End}{End}

 \DeclareMathOperator{\car}{char}

  \DeclareMathOperator{\e}{e}
  \DeclareMathOperator{\depth}{depth}

\def\urltilda{\kern -.15em\lower .7ex\hbox{\~{}}\kern .04em}
\def\urldot{\kern -.10em.\kern -.10em}\def\urlhttp{http\kern -.10em\lower -.1ex
\hbox{:}\kern -.12em\lower 0ex\hbox{/}\kern -.18em\lower 0ex\hbox{/}}

\selectfont

\begin{document}

\title[Brauer-Thrall for MCM Modules]{Brauer-Thrall Theory for\\ Maximal Cohen-Macaulay Modules}

\author{Graham J. Leuschke}
\address{Department of Mathematics\\
Syracuse University\\
Syracuse, NY 13244, USA}
\email{gjleusch@math.syr.edu}

\author{Roger Wiegand}
\address{Department of Mathematics \\ 
University of Nebraska--Lincoln \\ 
Lincoln, NE 68588, USA}
\email{rwiegand1@math.unl.edu}

\subjclass[2010]{%
Primary: 
13C14, 
16G50, 
16G60; 
Secondary: 
13C05, 
13H10, 
13J10, 
16G30
}


\keywords{Brauer-Thrall Conjectures, maximal Cohen-Macaulay modules, multiplicity}

\thanks{Research for this work was partially supported by NSF grant
  DMS-0902119 (GJL) and by a Simons Foundation Collaboration Grant (RW)}

\date{\today}


\begin{abstract} %
The Brauer-Thrall Conjectures, now theorems, were originally stated for finitely generated modules over a finite-dimensional $\sk$-algebra.  They say, roughly speaking, that infinite representation type implies the existence of lots of indecomposable modules of arbitrarily large $\sk$-dimension.  These conjectures have natural interpretations in the context of maximal Cohen-Macaulay modules over Cohen-Macaulay local rings.  This is a survey of progress on these transplanted conjectures.
\end{abstract}

\maketitle

\section{Introduction}

The Brauer-Thrall Conjectures first appeared in a 1957 paper by
Thrall's student J.\ P.\ Jans~\cite{Jans}.  They say, roughly
speaking, that if a finite-dimensional algebra $A$ over a field $\sk$
has infinite representation type, then $A$ has lots of big
indecomposable finitely generated modules.  Recall that $A$ has
\emph{finite representation type} provided there are only finitely
many indecomposable finitely generated $A$-modules up to isomorphism,
\emph{bounded representation type} provided there is a bound on the
$\sk$-dimensions of the indecomposable finitely generated $A$-modules,
and \emph{strongly unbounded representation type} provided there is an
infinite sequence $n_1<n_2<\cdots$ of positive integers such that
$A$ has, for each $i$, infinitely many non-isomorphic indecomposable
modules of $\sk$-dimension $n_i$.  Here are the conjectures:

\begin{conj}[First Brauer-Thrall Conjecture (BT1)]\label{conj:BT1} If
  $A$ has bounded representation type, then $A$ has finite
  representation type.
\end{conj}

\begin{conj}[Second Brauer-Thrall Conjecture (BT2)]\label{conj:BT2} If
  $A$ has unbounded representation type and $\sk$ is infinite, then $A$
  has strongly unbounded representation type.
\end{conj}

Under mild hypotheses, both conjectures are now theorems.  Ro\u\i
ter~\cite{Roiter:1968} verified (BT1), and Nazarova and Ro\u\i ter
\cite{Nazarova-Roiter:1973} proved (BT2) for perfect fields $\sk$.
See~\cite{Ringel:Report} or \cite{Gustafson:1982} for some history on
these results.

When we move from Artinian rings to  local rings $(R,\m,\sk)$ of
positive dimension, the first thing we need to do is to decide on the
right class of modules.  If $R$ is not a principal ideal ring,
constructions going back to Kronecker \cite{Kronecker:1874} and
Weierstra\ss\  \cite{Weierstrass:1868} show that $R$ has
indecomposable modules requiring arbitrarily many generators.
Moreover, if $\sk$ is infinite, then for every $n$ there are $|\sk|$
non-isomorphic indecomposable modules each of which requires exactly
$n$ generators.  (See \cite[Theorem 3.3 and Exercise 3.25]{BOOK}.)
Thus imposing finiteness or boundedness conditions on the class of
\emph{all} modules does not lead to anything interesting.  

Restricting to torsion-free modules yields a more robust theory, at
least in dimension one.  In the 1960's Jacobinski
\cite{Jacobinski:1967} and, independently, Drozd and Ro\u\i ter
\cite{Drozd-Roiter} studied orders in algebraic number fields and,
more generally, rings essentially module-finite over the ring of
integers, and classified the rings having only finitely many
indecomposable finitely generated torsion-free modules up to
isomorphism. 

 In dimensions greater than one, there are just too many
torsion-free modules. Indeed, Bass \cite{Bass:1962} proved in 1962
that every local domain of dimension two or more has indecomposable
finitely generated torsion-free modules of arbitrarily large rank.

The maximal Cohen-Macaulay (MCM) property, a higher-dimensional form
of torsion-freeness, turns out to give a fruitful class of modules to
study. The equality of a geometric invariant (dimension) with an
arithmetic one (depth) makes MCM modules easy to work with,
simultaneously ensuring that in some sense they faithfully reflect the
structure of the ring. For example, a Cohen-Macaulay local ring has no
non-free MCM modules if and only if it is a regular local ring, so the
rings that are the simplest homologically are also simple in this
sense.  Imposing finiteness or boundedness conditions on MCM modules
over a Cohen-Macaulay local ring leads to classes of rings that are
large enough to include interesting examples, but small enough to
study effectively. The seminal work of Herzog \cite{Herzog:1978},
Artin and Verdier \cite{Artin-Verdier:1985}, Auslander
\cite{Auslander:rationalsing}, and Buchweitz, Greuel, Kn\"orrer, and
Schreyer \cite{BGS, Knorrer} supports this assertion.  For example,
the main result of \cite{BGS, Knorrer} is that a complete
equicharacteristic hypersurface singularity over an algebraically
closed field of characteristics zero has only finitely many
indecomposable MCM modules up to isomorphism if and only if it is a
simple singularity in the sense of V.\ I.~Arnol$'$d, that is, one of
the (A$_n$), (D$_n$), (E$_6$), (E$_7$), or (E$_8$) hypersurface
singularities.

Next, we have to decide what invariant should be used to measure the size of a finitely generated module $M$.  Two obvious choices are $\mu_R(M)$, the minimal number of generators required for $M$, and $\e_R(M)$, the multiplicity of $M$.  We choose multiplicity. 

\begin{dfn}\label{dfn:bdd} Let $(R,\m,\sk)$ be a local ring. 
 \begin{enumerate}[label=(\roman{*})]
\item\label{item:finite} $R$ has \emph{finite} CM type provided $R$ has, up to isomorphism, only finitely many indecomposable MCM modules.
\item\label{item:bdd} $R$ has
 \emph{bounded} CM type provided there is a bound on the multiplicities of the indecomposable MCM $R$-modules.  
 \item\label{item:str-unbdd}$R$ has \emph{strongly unbounded} CM type provided there is an increasing sequence $n_1<n_2<\cdots$ of positive integers such that, for every $i$, there are infinitely many indecomposable MCM modules of multiplicity $n_i$. 
 \end{enumerate} 
\end{dfn}

Here, then, are the Brauer-Thrall Conjectures for MCM modules:

\begin{conj}[First Brauer-Thrall Conjecture for MCM modules (BTM1)]\label{conj:BTM1}  If a local  ring $(R,\m,\sk)$ has bounded CM type, then $R$ has finite CM type.
\end{conj}

\begin{conj}[Second Brauer-Thrall Conjecture for MCM modules (BTM2)]\label{conj:BTM2}   If a local  ring  $(R,\m,\sk)$ has  unbounded CM type and $\sk$ is infinite, then $R$ has strongly unbounded CM type.
\end{conj}

For MCM modules, multiplicity and number of generators enjoy a linear relationship:
\begin{equation}\label{eq:mu-mult}
\mu_R(M)\le \e_R(M) \le \e(R)\cdot\mu_R(M)\,,
\end{equation}
for every MCM $R$-module.  (See \cite[Corollary A.24]{BOOK} for a proof of the first inequality.)
It follows that we could replace multiplicity by number of generators in Definition~\ref{dfn:bdd} without changing the class of rings satisfying bounded (respectively, strongly unbounded) CM type.

In fact, Conjecture~\ref{conj:BTM1} is false, the designation ``conjecture'' being merely a convenient nod to history.  The first counterexample was given by Dieterich in 1980 \cite{Dieterich:1980}.  Let $\sk$ be a field of characteristic two, let $A = \sk[\![x]\!]$, and let $G$ be the two-element group.  Then $AG$  has bounded but infinite CM type.  Of course $AG$ is isomorphic to $\sk[\![x,y]\!]/(y^2)$, which, as we will see in the next section, has bounded but infinite CM type for \emph{any} field $\sk$.  

\begin{conv-not} Throughout, $R$ will be a local ring (always
  assumed to be commutative and Noetherian). The notation $(R,\m,\sk)$
  indicates that $\m$ is the maximal ideal of $R$ and that $\sk$ is
  the residue field $R/\m$.  All modules are assumed to be finitely
  generated.  The $\m$-adic completion of $R$ is $\widehat R$, and the
  integral closure of $R$ in its total quotient ring
  $K:=\{\text{non-zerodivisors}\}^{-1}R$ is $\overline R$.  A module
  $M$ is \emph{maximal Cohen-Macaulay} (abbreviated ``MCM'') provided
  $\depth(M) = \dim(R)$.  We will denote the multiplicity $\e(\m,M)$
  of the maximal ideal on $M$ simply by $\e_R(M)$, and we write
  $\e(R)$ instead of $\e_R(R)$.  (See \cite[Chapter 14]{Matsumura}.)
  The modifier ``Cohen-Macaulay'', when applied to the ring $R$, will
  often be abbreviated ``CM''.  Our standard reference for matters
  commutative-algebraic will be~\cite{Matsumura}, and for
  representation theory we refer to~\cite{BOOK} or~\cite{Yoshino:book}.
\end{conv-not}

\section{Dimension One}\label{sec:BTM-dim1} Before getting started, let's observe that both conjectures are true for local Artinian rings.  In this case all finitely generated modules are MCM modules.  If $(R,\m)$ is an Artinian principal ideal ring with $\m^t=0$, the indecomposable modules are $R/\m^i$, $1\le i \le t$.  We have already observed that if $R$ is not a principal ideal ring, then there exist, for each $n\ge1$, indecomposable modules requiring exactly $n$ generators, and, if $\sk$ is infinite, $|\sk|$ of them.

Now, on to dimension one!   We recall the characterization of one-dimensional rings of finite CM type:

\begin{thm}\label{thm:frt-dim1} Let $(R,\m,\sk)$ be a Cohen-Macaulay local ring of dimension one. Then $R$ has finite CM type if and only if
\begin{enumerate}[label=(\roman{*})]
\item $R$ is reduced,
\item $\mu_R(\overline R) \le 3$, and
\item $\frac{\m \overline R + R}{R}$ is cyclic as an $R$-module.
\end{enumerate}
\end{thm}

Items (i) and (ii)  are equivalent to the condition that $\widehat R$
is reduced and $\e(R) \le 3$.  Conditions (ii) and (iii)  are often
called the ``Drozd-Ro\u\i ter conditions'' \cite{CWW} to recognize the
1966 paper \cite{Drozd-Roiter} where they first appeared and were
shown to characterize the rings of finite CM type among local rings
essentially module-finite over $\ZZ$.  The work of Drozd and Ro\u\i ter was
clarified considerably in 1978 by Green and Reiner
\cite{Green-Reiner}, who used explicit matrix reductions to verify
finiteness of CM type in the presence of the Drozd-Ro\u\i ter
conditions.  In 1989 R.~Wiegand~\cite{Wiegand:1989} adapted constructions in \cite{Drozd-Roiter} to prove the ``only if'' direction in general.  A separable base-change argument in \cite{Wiegand:1989} and the matrix decompositions of Green and Reiner verified the ``if'' direction, except in the case of an imperfect residue field of characteristic two or three.  In \cite{Wiegand:1994} Wiegand took care of the case of characteristic three.  Finally, \c Cimen, in his Ph.D.\ dissertation \cite{Cimen:thesis}, completed the proof of Theorem~\ref{thm:frt-dim1} via difficult matrix reductions. (Cf.\ \cite{Cimen:paper}.)  

Although we won't say much about non-CM rings, we record the following result from \cite{Wiegand:1994}, which, together with Theorem~\ref{thm:frt-dim1}, characterizes the one-dimensional local rings with finite CM type:

\begin{thm}\label{non-CM}Let $(R,\m,\sk)$ be a one-dimensional local ring, not necessarily Cohen-Macaulay, and let $N$ be the nilradical of $R$.  Then $R$ has finite CM type if and only if
\begin{enumerate}[label=(\roman{*})]
\item $R/N$ (which \emph{is} CM) has finite CM type, and
\item $N\cap \m ^n = 0$ for some positive integer $n$.
\end{enumerate}
\end{thm}

The proof of the ``only if'' direction in Theorem~\ref{thm:frt-dim1} (necessity of the Drozd-Ro\u\i ter conditions) in \cite{Wiegand:1989} actually shows more and confirms BTM1 in the analytically unramified case.  We will say that a finitely generated module $M$ over a CM local ring $R$ has \emph{constant rank} $n$ provided $K\otimes_RM\cong K^{(n)}$, where $K$ is the total quotient ring.  Equivalently, $M_\p$ is a free $R_\p$-module of rank $n$ for every minimal prime ideal $\p$ of $R$.   In this case $\e(M) = n\e(R)$.

\begin{thm}[BTM1 when $\widehat R$ is reduced, \cite{Wiegand:1989}]\label{thm:BTM1-red} Let $(R,
\m,\sk)$ be a one-dimensional local ring with reduced completion.  If $R$ has infinite CM type then for every $n$ there is an indecomposable MCM $R$-module of constant rank $n$.  In particular, $R$ has unbounded CM type.
\end{thm}

We have already seen that BTM1 can fail if there are nilpotents.  We
showed in 2005 \cite[Theorem 2.4]{Leuschke-Wiegand:bcmt} that there
are essentially only three counterexamples to BTM1 in dimension one.
Recall that the complete (A$_\infty$) and (D$_\infty$) curve
singularities are, respectively, the rings $\sk[\![x,y]\!]/(y^2)$ and
$\sk[\![x,y]\!]/(xy^2)$.  They arise as the respective limits of the (A$
_n$) singularities $\sk[\![x,y]\!]/(y^2+x^{n+1})$ and the (D$_n$)
singularities $\sk[\![x,y]\!]/(xy^2+x^{n-1})$ as $n\to\infty$.

\begin{thm}[Failure of BTM1, \cite{Leuschke-Wiegand:bcmt}]\label{thm:BTM1-fails}
Let $(R,\m,\sk)$ be an equicharacteristic, one-dimensional, Cohen-Macaulay local ring, with $\sk$ infinite.  Then $R$ has bounded but infinite CM type if and only if the completion $\widehat R$ is isomorphic to one of the following:
\begin{enumerate}[label=(\roman{*})]
\item $\sk[\![x,y]\!]/(y^2)$, the (A$_\infty$) singularity;
\item $T:= \sk[\![x,y]\!]/(xy^2)$, the (D$_\infty$) singularity;
\item $E :=\End_T(\m_T)$, the endomorphism ring of the maximal ideal of $T$.
\end{enumerate}
The ring $E$ has a presentation $E\cong \sk[\![X,Y,Z]\!]/(XY,YZ,Z^2)$.
\end{thm}

The assumption that $\sk$ be infinite is annoying.  It's tempting to
try to eliminate this assumption via the flat local homomorphism $R
\to S:=R[z]_{\m R[z]}$, where $z$ is an indeterminate.  The problem
would be to show that if $R$ has unbounded CM type then so has $S$.
While it is rather easy to show that finite CM type descends along
flat local homomorphisms (as long as the closed fiber is CM)
\cite[Theorem 1.6]{Wiegand:1998}, it's not known (at least to us)
whether an analogous result holds for descent of \emph{bounded} CM
type.  In fact, it is not even known, in higher dimensions, whether
bounded CM type descends from the completion.  Such descent was a
crucial part of the proof of Theorem~\ref{thm:BTM1-fails}, but the
proof of descent was based not on abstract considerations, but on the
precise presentations, in \cite{BGS}, of  the indecomposable $\widehat
R$-modules in each of the three cases.  Using these presentations, we
were able to say exactly which MCM $\widehat R$-modules are extended
from $R$-modules, and thereby deduce that $R$ itself has bounded CM
type.  Part of the difficulty in proving a general statement of this
form is that there may be no uniform bound on the number of indecomposable MCM
$\widehat R$-modules required to decompose the completion of an
indecomposable MCM $R$-module (see~\cite[Example 17.11]{BOOK}).  

\medskip

At this point we have shown, for CM local rings of dimension one, that BTM1 holds in the analytically unramified case but fails (just a little bit) in general.  We turn now to BTM2 for CM local rings of dimension one.  

\begin{thm}\label{BTM2-dim1}Let $(R,\m,\sk)$ be a one-dimensional local Cohen-Macaulay ring with unbounded CM type, and with $\sk$ infinite.  Assume either
\begin{enumerate}[label=(\roman{*})]
\item $\widehat R$ is reduced, or
\item $R$ contains a field.
\end{enumerate}
Then, for each positive integer $n$, $R$ has $|\sk|$ pairwise non-isomorphic indecomposable MCM modules of constant rank $n$.  In particular, BTM2 holds for one-dimensional CM local rings that satisfy either (i) or (ii).
\end{thm}

Karr and R.~Wiegand \cite[Theorem 1.4]{Karr-Wiegand:BT2} proved this in the analytically unramified case (i).  Later Leuschke and Wiegand modified that proof, using ideas from \cite{Leuschke-Wiegand:hyperbrt} and \cite{Leuschke-Wiegand:bcmt}, to prove the result in the equicharacteristic case (ii). See \cite[Theorem 17.10]{BOOK}.  The rest of this section is devoted to a sketch of the main ideas of the proof of Theorem~\ref{BTM2-dim1}.

Assume, for the rest of this section, that $(R,\m,\sk)$ is a one-dimensional CM local ring satisfying the hypotheses of Theorem~\ref{BTM2-dim1}.  In particular,  $\sk$ is  infinite and $R$ has
unbounded CM type. The first step, proved by Bass \cite{Bass:ubiquity} in the analytically unramified case, appears as Theorem 2.1 of \cite{Leuschke-Wiegand:hyperbrt}:

\begin{lem}\label{lem:mult2}Suppose $\e(R) \le 2$.  Then every indecomposable MCM $R$-module is isomorphic to an ideal of $R$ and hence has multiplicity at most two.
\end{lem}

Thus we may assume that $\e(R)\ge 3$, and in this case $R$ has a finite birational extension $S$ (an intermediate ring between $R$ and its total quotient ring $K$ such that $S$ is finitely generated as an $R$-module) with $\mu_R(S) = \e(R)$.  Although we will need to choose $S$ with some care, we note here that  $S := \bigcup_{n\ge 1}\End_R(\m^n)$ has the right number of generators. (See \cite[Lemma 2.6]{Leuschke-Wiegand:hyperbrt}.)  In the analytically unramified case, one typically takes $S = \overline R$.  (Notice that none of this works if $R$ is not CM, since  $R=K$ in that case!)  Let $\ff$ be the conductor, that is, the largest ideal of $S$ that is contained in $R$.  Putting $A=R/\ff$, $B=S/\ff$, and $D=B/\m B$, we obtain a commutative diagram
\begin{equation}\label{cond-sq}
\begin{gathered}
\xymatrix{
R \ar[r] \ar[d] & S \ar[d] \\
A \ar[r] \ar[d] & B \ar[d] \\
\sk \ar[r] & D
}
\end{gathered}
\end{equation}
in which the top square is a pullback and $D$ is a $\sk$-algebra of dimension $\e(R)$.

Now let $n$ be a fixed positive integer, and let $t\in \sk$.  We wish to build a family, parametrized by $t$, of indecomposable MCM $R$-modules of constant rank $n$.  The following  construction \cite[Construction 2.5]{Wiegand:1989}, \cite[Construction 3.13]{BOOK} is based on work of Drozd and Ro\u\i ter \cite{Drozd-Roiter}.    Let $I$ be the $n\times n$ identity matrix and $H$ be the nilpotent $n\times n$ Jordan block with $1$'s on the superdiagonal and $0$'s elsewhere.  Let  $\alpha$ and  $\beta$ be elements of  $D$ such that $\{1,\alpha,\beta\}$ is linearly independent over $\sk$.  (Eventually we will have to impose additional restrictions on $\alpha$ and $\beta$.) Let $V_t$ be the $\sk$-subspace of $D^{(n)}$ spanned by the columns of the $n \times 2n$ matrix
\begin{equation}\label{eq:Psi-def}
\Psi_t:= \left[I \quad \alpha I + \beta(tI+H)\right]\,. 
\end{equation}
Let $\pi\colon S^{(n)} \twoheadrightarrow D^{(n)}$ be the canonical surjection, and define $M_t$ by the following pullback diagram.
\begin{equation}\label{CD:M_t}
\begin{gathered}
\xymatrix{
M_t \ar[r] \ar[d] & S^{(n)} \ar[d]^\pi \\
V_t \ar@{^{(}->}[r] & D^{(n)}
}
\end{gathered}
\end{equation}

Then $M_t$ is an MCM $R$-module of constant rank $n$, and it is indecomposable if the \emph{pair} $V_t\subseteq D^{(n)}$ is indecomposable in the following sense:  There is no idempotent endomorphism $\varepsilon$ of $D^{(n)}$, other than $0$ and the identity, such that $\varepsilon(V_t)\subseteq V_t$.  Moreover, if $t,u\in \sk$ and $M_t\cong M_u$, then the pairs $(V_t\subseteq D^{(n)})$ and $(V_u\subseteq D^{(n)})$ are isomorphic, in the sense that there is an automorphism $\varphi$ of $D^{(n)}$ such that $\varphi(V_t)\subseteq V_u$.  Our goal, then, is to choose $\alpha$ and $\beta$ so that we get $|\sk|$ non-isomorphic indecomposable pairs $(V_t\subseteq D^{(n)})$.  

Suppose first that $\e(R) = 3$.  We need to choose a finite birational extension $R\subset S$ such that
\begin{equation}\label{eq:mult-3}
\mu_R(S)=3 \quad \text{and} \quad \mu_R\left(\frac{\m S + R}{R}\right) \ge 2\,.
\end{equation}

If $R$ is analytically unramified, the assumption that $R$ has
unbounded (hence infinite) CM type implies failure of the second
Drozd-Ro\u\i ter condition (iii) in Theorem~\ref{thm:frt-dim1}, and we
can take $S = \overline R$.  If $R$ is analytically ramified but
contains a field, the fact that $\widehat R$ is \emph{not} one of the
three exceptional rings of Theorem~\ref{thm:BTM1-fails} leads, after
substantial computation, to the right choice for $S$.  (See the
proof of \cite[Theorem 1.5]{Leuschke-Wiegand:bcmt} or the proofs of Theorems 17.6 and 17.9 in \cite{BOOK}.)

Now, with our carefully chosen birational extension $R\to S$, we have
\begin{equation}\label{eq:mult-3-art}
\dim_\sk(B/\m B) = 3 \quad \text{and} \quad \dim_\sk\left(\frac{\m B +A}{\m^2 B+A}\right) \ge 2\,,
\end{equation}
for the Artinian rings $A$ and $B$ in the diagram~\eqref{cond-sq}.  Put $C = \m B+A$, and choose elements $x,y\in \m B$ so that their images in $\frac{\m B +A}{\m^2 B+A}$ are linearly independent.  Since $C/\m C$ maps onto $\frac{\m B +A}{\m^2 B+A}$, the images $\alpha$ and $\beta$ of $x$ and $y$ in $C/\m C$ are linearly independent.  By \cite[Lemmas 3.10 and 3.11]{BOOK} it suffices to build the requisite pairs $(V_t\subseteq (C/\m C)^{(n)})$, since these will yield, via extension, non-isomorphic indecomposable pairs $(V_t\subseteq D^{(n)})$.  Moreover, with this choice of $\alpha$ and $\beta$, we have the relations
\begin{equation}\label{eq:very-short}
\alpha^2 = \alpha\beta = \beta^2 = 0\,.
\end{equation}

Returning to the general case $\e(R)\ge 3$, we may assume that either $\dim_\sk(D) \ge 4$ or else $D$ contains elements $\alpha$ and $\beta$ satisfying \eqref{eq:very-short}.  In order to show that there are enough values of $t$ that produce non-isomorphic indecomposable pairs $(V_t\subseteq D^{(n)})$, we let $t$ and $u$ be elements of $\sk$, not necessarily distinct,  and suppose that $\varphi$ is a $\sk$-endomorphism of $D^{(n)}$ that carries $V_t$ into $V_u$. We regard $\varphi$ as an $n\times n$ matrix with entries in $D$.  Recalling that $V_t$ is the column space of the matrix $\Psi_t$ in \eqref{eq:Psi-def}, we see that the condition $\varphi V_t\subseteq V_u$ yields a $2n\times 2n$ matrix $\theta$ over $\sk$ satisfying the equation
\begin{equation}\label{eq:commute}
\varphi\Psi_t =\Psi_u\theta\,.
\end{equation}
\noindent Write $\theta = \left[
  \begin{smallmatrix} E&F\\
    P&Q
  \end{smallmatrix}\right]$, where $E$, $F$, $P$, and  $Q$
are $n\times n$ blocks.  Then~\eqref{eq:commute} gives the following two
equations:
\begin{equation}\label{eq:blarg}
  \begin{aligned}
    \varphi &=E+\alpha P+\beta(u I+H)P\\
    \alpha\varphi + \beta\varphi(t I+H) &=F+\alpha Q+\beta(u I+H)Q\,.
\end{aligned}
\end{equation}
Substituting the first equation into the second and combining terms,
we get the following equation:
\begin{multline}\label{eq:mess}
  -F + \alpha(E-Q) + \beta(t E-u Q + E H - HQ) +
  (\alpha+t\beta)(\alpha+u\beta)P \\
  + \alpha\beta(HP + P H) + \beta^2(H P H + t H P+u P H) = 0\,.
\end{multline}

Suppose there exist elements $\alpha$ and $\beta$ satisfying  Equation~\eqref{eq:very-short}.  With this choice of $\alpha$ and $\beta$, \eqref{eq:mess} collapses:
\begin{equation}\label{eq:gleep}
-F + \alpha(E-Q) + \beta(t E-u Q + E H - HQ)  =0\,.
\end{equation}
 Since all capital letters in \eqref{eq:gleep} represent matrices over $\sk$, and since $\{1,\alpha,\beta\}$ is linearly independent over $\sk$, we get the equations
$$
F=0\,, \qquad E=Q\,, \qquad \text{and} \qquad (t-u) E + EH-HE = 0\,.
$$
  After a bit of fiddling (see \cite[Case 3.14]{BOOK} for the details) we reach two conclusions:  
\begin{enumerate}[label=(\roman{*})]
\item  If $\varphi$ is invertible, then $t = u$.  Thus the modules are pairwise non-isomorphic.
\item If $t=u$ and $\varphi$ is idempotent, then $\varphi$ is either $0$ or $I$.  Thus all of the modules are indecomposable.   
\end{enumerate}
The key issue in these computations is that the matrix $H$ is non-derogatory, so that its commutator in the full matrix ring is just the local ring $\sk[H]\cong \sk[X]/(X^n)$.   

We may therefore assume that $\dim_\sk(D) \ge 4$.   With a little luck, the algebra $D$ might contain an element $\alpha$ that does \emph{not} satisfy a non-trivial quadratic relation over $\sk$.  In this case, we choose any element $\beta\in D$ so that $\{1,\alpha,\alpha^2,\beta\}$ is linearly independent, and we set
\[
G = \{t\in \sk\mid \{1,\alpha,\beta,(\alpha+t\beta)^2\} \text{ is  linearly independent}\}\,.
\]
This set is non-empty and Zariski-open, hence cofinite in $\sk$.  For $t\in G$, put
\[
G_t = \{u\in G \mid \{1,\alpha,\beta,(\alpha+t\beta)(\alpha+u\beta)\} \text{ is  linearly independent}\}\,.
\]
Then $G_t$ is cofinite in $G$ for each $t\in G$.  Moreover, one can check the following, using \eqref{eq:mess}: 
\begin{enumerate}[label=(\roman{*})]
\item If $t$ and $u$ are distinct elements of $G$ with $u\in G_t$, then $\varphi$ is not an isomorphism.
\item If $t=u\in G$ and $\varphi$ is idempotent, then $\varphi$ is either $0$ or $I$.
\end{enumerate}
\noindent The desired conclusions follow easily.  (See \cite[Case 3.16]{BOOK} for the details.)

The remainder of the proof \cite[(3.17)--(3.21)]{BOOK} is a careful analysis of the $\sk$-algebras $D$ in which every element is quadratic over $\sk$.  (The fact that $\sk$ is infinite obviates consideration of the last case \cite[Case 3.22]{BOOK}, where our  construction does not work and  Dade's construction \cite{Dade:1963} is used to produce \emph{one} indecomposable of rank $n$.)

\medskip

In studying direct-sum decompositions over one-dimensional local
rings, it is important to know about indecomposable MCM modules of
non-constant rank. (See \cite{Wiegand:1991}, where Wiegand determined exactly how badly Krull-Remak-Schmidt uniqueness can fail.)  If $(R,\m,\sk)$ is a one-dimensional, analytically unramified local ring with minimal prime ideals $\p_1,\dots,\p_s$, we define the \emph{rank} of a module to be the $s$-tuple
$(r_1,\dots,r_s)$, where $r_i$ is the dimension of $(M_{\p_i})$ as a vector space over the
field  $R_{\p_i}$.  Crabbe and Saccon \cite{Crabbe-Saccon} have recently
proved the following:

\begin{thm}\label{thm:Crabbe-Saccon} Let $(R,\m,\sk)$ be an analytically unramified local ring of dimension one, with minimal prime ideals $\p_1,\dots,\p_s$.  Assume that $R/\p_1$ has infinite CM type.  Let $\underline r := (r_1,\dots,r_s)$ be an arbitrary $s$-tuple of non-negative integers with $r_1\ge r_i$ for each $i$, and with $r_1>0$.  Then there is an indecomposable MCM $R$-module with $\rank(M) = \underline r$, and $|\sk|$ non-isomorphic ones if $\sk$ is infinite.
\end{thm}

\section{Brauer-Thrall I for Hypersurfaces}

In Theorem~\ref{thm:BTM1-fails} we saw that there are just two plane curve singularities that contradict BTM1.  Here we promote this result to higher-dimensional hypersurfaces,  with the following theorem from \cite{Leuschke-Wiegand:hyperbrt} (cf.\ \cite[Theorem 17.5]{BOOK}): 
\begin{thm}\label{thm:hyperbrt} Let $\sk$ be an algebraically closed
  field of characteristic different from two, and let $R =
  \sk[\![x_0,\dots,x_d]\!]/(f)$, where $f$ is a non-zero element of
  $(x_0,\dots,x_d)$ and $d\ge 2$. Then $R$ has bounded but infinite
  CM type if and only if  $R\cong \sk[\![x_0,\dots,
  x_d]\!]/(g+x_2^2+\cdots+x_d^2)$, where $g$ is a polynomial in
  $\sk[x_0,x_1]$ defining either an (A$_\infty$) or (D$_\infty$) curve singularity.
\end{thm}

This theorem and its proof are modeled on the beautiful result of Buchweitz, Greuel, Kn\"orrer, and Schreyer, where ``bounded but infinite'' is replaced by ``finite'', and the singularities in the conclusion are the \emph{simple} or \emph{ADE} singularities, \cite[\S4.3]{BOOK}.  

The ``if'' direction of Theorem~\ref{thm:hyperbrt} hinges on the following result (see \cite[Theorem 17.2]{BOOK}):

\begin{lem}[Kn\"orrer, \cite{Knorrer}]\label{lem:Knorrer} Let $\sk$ be a field, and put $S=\sk[\![x_0,\dots,x_d]\!]$.  Let $f$ be a non-zero non-unit of $S$, $R=S/(f)$, and $R^\# = S[\![z]\!]/(f+z^2)$.  
\begin{enumerate}[label=(\roman{*})]
\item If $R^\#$ has finite (respectively, bounded) CM type, so has $R$.
\item Assume $R$ has finite (respectively, bounded) CM type and $\car(\sk)\ne2$.  Then $R^\#$ has finite (respectively, bounded) CM type.  More precisely, if $\mu_R(M) \le B$ for every indecomposable MCM $R$-module $M$, then $\mu_{R^\#}(N) \le 2B$ for every indecomposable MCM $R^\#$-module $N$.
\end{enumerate}
\end{lem}

For the ``only if'' direction, we need Lemma~\ref{lem:Knorrer} and the following result due to Kawasaki \cite[Theorem 4.1]{Kawasaki:1996}:

\begin{lem}\label{Kawasaki} Let $(R,\m)$ be a $d$-dimensional abstract hypersurface (a local ring whose completion $\widehat R$ has the form $S/(f)$, where $(S,\n)$ is a regular local ring and $f\in \n$).  Let $n$ be any positive integer, and let $M$ be the $(d+1)^{\text{st}}$ syzygy of $R/\m^n$.  If $\e(R) > 2$, then $M$ is an indecomposable MCM $R$-module, and $\mu_R(M) \ge \binom{d+n-1}{d-1}$.  In particular, if $d\ge 2$ then $R$ has unbounded CM type.
\end{lem} 

If, now, $d\ge 2$ and $R$ (as in Theorem~\ref{thm:hyperbrt}) has bounded but infinite CM type, then $\e(R) \le 2$.  Using the Weierstra\ss\ Preparation Theorem and a change of variables, we can put $f$ into the form $g+x_d^2$, with $g\in \sk[\![x_0,\dots,x_{d-1}]\!]$.  Then $\sk[\![x_0,\dots,x_{d-1}]\!]/(g)$ has bounded but infinite CM type, by Lemma~\ref{lem:Knorrer}.  We repeat this process till we get down to dimension one, and then we invoke Theorem~\ref{thm:BTM1-fails}.

\section{Brauer-Thrall I for Excellent Isolated Singularities}

The starting point here is the Harada-Sai Lemma \cite[Lemmas 11 and 12]{Harada-Sai}, sharpened by Eisenbud and de la Pe\~na in 1998 \cite{Eisenbud-delaPena:1998}.  By a \emph{Harada-Sai sequence} we mean a sequence
\[
M_1\overset{f_1}{\to}M_2\overset{f_2}{\to}\cdots\overset{f_{s-1}}{\to}M_s
\]
of $R$-homomorphisms satisfying
\begin{enumerate}[label=(\roman{*})]
\item each $M_i$ is indecomposable of finite length;
\item no $f_i$ is an isomorphism; and
\item the composition $f_{s-1}f_{s-2}\cdots f_1$ is non-zero.
\end{enumerate}

\begin{lem}[Harada-Sai Lemma] With the notation above, suppose $\ell_R(M_i) \le b$ for each $i$.  Then $s\le 2^b - 1$.
\end{lem}

In fact, Eisenbud and de la Pe\~na \cite{Eisenbud-delaPena:1998} characterized the integer sequences  that can occur in the form $(\ell_R(M_1),\dots, \ell_R(M_s))$ for some Harada-Sai sequence over some $R$.  In order to apply Harada-Sai to MCM modules, we need to reduce modulo a suitable system of parameters to get down to the Artinian case.  Of course, an arbitrary system of parameters won't work, since indecomposability and non-isomorphism won't be preserved.  What we need is a \emph{faithful system of parameters}, that is, a system of parameters
$\underline x = x_1,\dots, x_d$ such that $\underline x\Ext^1_R(M,N) = 0$ for every MCM $R$-module $M$ and every finitely generated $R$-module $N$.  Here are some useful properties of faithful systems of parameters (where we write $\underline x^2$ for the system of parameters ($x_1^2,\dots,x_d^2$)):

\begin{lem}\label{lem:faithful-sop} Let $\underline x$ be a faithful system of parameters for a CM local ring $R$.
\begin{enumerate}[label=(\roman{*})]
\item Let $M$ and $N$  be MCM $R$-modules, and suppose
  $\varphi\colon M/\underline x^2M \to N/\underline x^2N$ is an isomorphism.
  There is an isomorphism $\tilde\varphi\colon M\to N$ such that $\tilde\varphi\otimes_R(R/(\underline x) )= \varphi\otimes_R(R/(\underline x))$.
\item Let $s:\ 0\to N \to E \to M \to 0$ be a short exact of MCM $R$-modules.  Then $s$ splits if and only if $s\otimes_R(R/(\underline x^2))$ splits.
\item Assume $R$ is Henselian, and let $M$ be an indecomposable MCM $R$-module.  Then $M/\underline x^2 M$ is indecomposable.
\end{enumerate}
\end{lem}

Using these properties, one obtains the Harada-Sai Lemma for MCM modules \cite[Theorem 15.19]{BOOK} 

\begin{lem}\label{Harada-Sai-MCM} Let $(R,\m,\sk)$ be a CM, Henselian local ring
and $\underline x$ a faithful system of parameters.  Let
\[
M_1\overset{f_1}{\to}M_2\overset{f_2}{\to}\cdots\overset{f_{s-1}}{\to}M_s
\]
be a sequence of $R$-homomorphisms, with each $M_i$ indecomposable and
MCM\@.  Assume that 
\[
(f_{s-1}f_{s-2}\cdots f_1)\otimes_R(R/(\underline x^2)) \ne 0\,.
\]
  If 
$\ell_R(M_i/\underline x M_i)\le b$ for all $i$, then $s\le 2^b-1$.
\end{lem}
Suppose, instead, that we have a bound, say, $B$, on the multiplicities $\e(M_i)$.  Choosing $t$ such that $\m^t\subseteq (\underline x^2)$, we get a bound $b:=Bt^d$ on the lengths of the modules $M_i/\underline x^2M_i$. A  walk around the AR quiver of $R$ then proves BTM1.  (See \cite[Chapter 6]{Yoshino:book} or \cite[\S15.3]{BOOK}.).  Of course, none of this does any good unless the ring $R$ \emph{has} a faithful system of parameters.  The big theorem here is due to Yoshino \cite{Yoshino:1987} (cf.\ \cite[Theorem 15.8]{BOOK}):

\begin{thm}\label{faithful-exist}Let $(R,\m,\sk)$ be a complete CM local ring containing a field.  Assume $\sk$ is perfect and that $R$ has an isolated singularity.  Then $R$ has a faithful system of parameters.
\end{thm}

Putting all of this stuff together, we obtain the following theorem, proved independently by Dieterich \cite{Dieterich:1987} and Yoshino \cite{Yoshino:1987}:

\begin{thm}\label{complete-BTM1} Let $(R,\m,\sk)$ be a complete, equicharacteristic local ring with perfect residue field $\sk$.  Then $R$ has finite CM type if and only if
\begin{enumerate}[label=(\roman{*})]
\item $R$ has bounded CM type, and
\item $R$ has an isolated singularity.
\end{enumerate}
\end{thm}

The main thrust is the ``if'' direction, the converse being a consequence of Auslander's famous theorem \cite{Auslander:isolsing} that complete CM rings with finite CM type must be isolated singularities.

In 2005, Leuschke and Wiegand used ascent and descent techniques to prove the following generalization \cite[Theorem 3.4]{Leuschke-Wiegand:bcmt}:

\begin{thm}\label{excellent-BTM1} Let $(R,\m,\sk)$ be an excellent, equicharacteristic local ring with perfect residue field $\sk$.  Then $R$ has finite CM type if and only if 
\begin{enumerate}[label=(\roman{*})]
\item $R$ has bounded CM type, and
\item $R$ has an isolated singularity.
\end{enumerate}
\end{thm}

This time, for the ``only if'' direction, one needs the Huneke-Leuschke version \cite{Huneke-Leuschke:2002} of Auslander's theorem, stating that \emph{every} CM ring of finite CM type has an isolated singularity.

Without the word ``excellent'', Theorem~\ref{excellent-BTM1} would be false.  For example, the ring $\mathbb C[\![x,y]\!]/(y^2)$ is the completion of an integral domain $(R,\m)$, by Lech's Theorem \cite{Lech}.  Theorem~\ref{thm:BTM1-fails} implies that $R$ has bounded but infinite CM type, and of course $R$ has an isolated singularity.

\section{Brauer-Thrall II} In Section~\ref{sec:BTM-dim1} we proved a strong form of BTM2 for one-dimensional CM local rings, assuming only that the ring is either analytically unramified or equicharacteristic.  In higher dimensions, no such general results are known.  One problem, already mentioned, is that there is no general result showing descent of bounded CM type along flat local homomorphisms.   Typically, one restricts to complete (or at least excellent Henselian) isolated singularities with algebraically closed residue field, in order to make use of the Auslander-Reiten quiver.

The following result was proved by Dieterich \cite[Theorem 20]{Dieterich:1987} in 1987, for characteristics different from two.  The case $\car(\sk) = 2$ was proved by  Popescu and Roczen \cite{Popescu-Roczen:1991} in 1991.

\begin{thm}\label{thm:BTM2-hyper}Let $R = \sk[\![x_0,\dots,x_d]\!]/(f)$ be a hypersurface isolated singularity, with $\sk$ algebraically closed.  If $R$ has infinite CM type, then $R$ has strictly unbounded CM type.
\end{thm}

Using Elkik's theorem \cite{Elkik} on modules extended from the Henselization, one can generalize this result to excellent Henselian rings (cf.\ \cite{Popescu-Roczen:1990}):

\begin{cor}\label{cor:BTM2-hensel}Let $(R,\m,\sk)$ be an excellent, equicharacteristic, Henselian local ring whose completion is a hypersurface.  Assume that $R$ has an isolated singularity and that $\sk$ is algebraically closed.  If $R$ has infinite CM type, then $R$ has strictly unbounded CM type.  (In particular, both BTM1 and BTM2 hold for these rings.)
\end{cor}

Excellence guarantees that the completion $\widehat R$ is an isolated singularity too.  (In fact, all one needs is that the inclusion $R\to \widehat R$ be a regular homomorphism (see \cite[Proposition 10.9]{BOOK}).)  If $N$ is an MCM $\widehat R$-module, then $N$ is free on the punctured spectrum of $\widehat R$ and hence, by \cite{Elkik}, is extended from an $R$-module.  This means that the map $M\mapsto \widehat M$, from MCM $R$-modules to MCM $\widehat R$-modules, is bijective on isomorphism classes.  Since $\e_R(M) = \e_{\widehat R}(\widehat M)$, the corollary follows from the theorem.

\medskip

The main thing we want to talk about in this section is Smal\o's remarkable result \cite{Smalo:1980} that produces, from an infinite family of indecomposable MCM modules of \emph{one fixed} multiplicity $n$, an integer $n' > n$ and an infinite family of indecomposable MCM modules of multiplicity $n'$.   In principle, this ought to make proofs of BTM2 lots easier.  We will give two such applications and also point out some limitations to this approach.  Here is Smal\o's theorem, proved in 1980 for Artin algebras:

\begin{thm}\label{Smalo} Let $(R,\m,\sk)$ be a complete CM isolated singularity, with $\sk$ algebraically closed.  Suppose $\{M_i\}_{i\in I}$ is an infinite family of pairwise non-isomorphic indecomposable MCM $R$-modules, all of the same multiplicity $n$.  There exist an integer $n'> n$, a subset $J$ of $I$ with $|J| = |I|$, and a family $\{N_j\}_{j\in J}$ of pairwise non-isomorphic indecomposable MCM $R$-modules, each of multiplicity $n'$.
\end{thm}

The basic ideas of Smal\o's proof survive transplantation to the MCM context remarkably well.  The proof uses the Harada-Sai Lemma~\ref{Harada-Sai-MCM} as well as a couple of lemmas that control multiplicity as one wanders around the AR quiver.  One  of these~\cite[Lemma 4.2.7]{Avramov:6lectures} bounds the growth of the Betti numbers $\beta_i(M)$ of a MCM module $M$ over a CM local ring of  multiplicity $e$:  $\beta_{i+1} \le (e-1) \beta_i$ for all $i$.   Another gives a linear bound between the multiplicities of the source and target of an irreducible homomorphism:  With $R$ as in the theorem, there is a positive constant $c$  such that $\e_R(M)\le c\e_R(N)\le c^2 \e_R(M)$ whenever $M\to N$ is an irreducible homomorphism of indecomposable MCM $R$-modules.  We refer the reader to~\cite[Section 15.4]{BOOK} for the details.

Here is an obvious corollary of Smal\o's theorem:

\begin{cor}\label{cor:uncountable} Let $(R,\m,\sk)$ be a complete CM isolated singularity, with $\sk$ algebraically closed.  If $R$ has uncountable CM type, then there is an sequence $n_1<n_2<n_3<\dots$ of positive integers such that $R$ has, for each $i$, uncountably many non-isomorphic indecomposable MCM modules of multiplicity $n_i$.
\end{cor}

As another application, one can give a proof of BTM2 in dimension one that is much less computational than the one given in Section~\ref{sec:BTM-dim1}, at least in an important  special case.  Suppose that $(R,\m,\sk)$ is a complete, reduced local ring of dimension one, and assume $R$ has infinite CM type.  Then the Drozd-Ro\u\i ter conditions ((ii) and (iii) of Theorem~\ref{thm:frt-dim1}) fail.  It is now a comparatively simple matter (see \cite[\S4]{Wiegand:1989}) to show that $R$ has an infinite family of pairwise non-isomorphic ideals.  We decompose each of these ideals into indecomposable summands, noting that $\e(R)$ bounds the number of summands of each ideal.  This yields infinitely many pairwise non-isomorphic indecomposable MCM modules, each with multiplicity  bounded by $\e(R)$, and hence an infinite subfamily consisting of modules of fixed multiplicity.  Now Smal\o's theorem shows that BTM2 holds for these rings.  

Don't be misled by this example.  In higher dimensions there is no hope of starting the inductive hypothesis with modules of rank one, in view of the following theorem due to Bruns \cite[Corollary 2]{Bruns:1981}:

\begin{thm}\label{thm:Bruns}  Let $A$ be any commutative Noetherian
  ring and $M$  a finitely generated $R$-module of constant rank $r$.
  Let $N$ be a second syzygy of $M$, and let $s$ be the (constant) rank of $N$.  If $M$ is not free, then the codimension of its non-free locus is at most $r+s+1$.  
\end{thm}

\begin{cor}\label{cor:Bruns-hyper} Let $(R,\m)$ be a $d$-dimensional isolated singularity whose completion is a hypersurface.  Let $M$ be a non-free MCM $R$-module of constant rank $r$.  Then $r\ge \frac{1}{2}(d-1)$.
\end{cor}

This bound is probably much too low.  In fact, Buchweitz, Greuel, and Schreyer \cite{BGS} conjecture that $r\ge 2^{d-1}$.  Nonetheless, the bound given in the corollary rules out MCM ideals once the dimension exceeds three. 

\section{Open Questions}\label{sec:questions}

Here we list a few open questions, some of which have already been mentioned at least implicitly.

\begin{ques}\label{ques:BRT2} Are there \emph{any} counterexamples to
  BTM2\@: Of course this is the same as asking whether BTM2 is true,
  but let's not even assume that $(R,\m,\sk)$ is CM\@.  What if
  $\dim(R) = 1$?  What if $\dim(R) = 1$ and $R$ is not CM\@?  The list
  goes on\dots. \end{ques}

\begin{ques}\label{ques:infinite-BTM1} Can one delete the assumption, in Theorem~\ref{thm:BTM1-fails}, that $\sk$ be infinite?
\end{ques}

\begin{ques}\label{ques:BCMT-descends} If $(R,\m)$ is a local CM ring whose completion $\widehat R$ has bounded CM type, must $R$ have bounded CM type?  More generally, let $R\to S$ be a flat local homomorphism with CM closed fiber.  If $S$ has bounded CM type, must $R$ have bounded CM type?
\end{ques}

\begin{ques}\label{ques:imperfect-BTM1} Can one delete the assumption, in Theorem~\ref{excellent-BTM1}, that $\sk$ be perfect?
\end{ques} 

\begin{ques} Can we improve Corollary~\ref{cor:Bruns-hyper}, getting better lower bounds for the rank, or multiplicity, of a non-free MCM module?
\end{ques}


\newcommand{\arxiv}[2][AC]{\mbox{\href{http://arxiv.org/abs/#2}{\textsf{arXiv:#2}}}}
\newcommand{\oldarxiv}[2][AC]{\mbox{\href{http://arxiv.org/abs/math/#2}{\textsf{arXiv:math/#2[math.#1]}}}}
\renewcommand{\MR}[1]{%
  {\href{http://www.ams.org/mathscinet-getitem?mr=#1}{MR #1}.}}
\providecommand{\bysame}{\leavevmode\hbox to3em{\hrulefill}\thinspace}
\newcommand{\arXiv}[1]{%
  \relax\ifhmode\unskip\space\fi\href{http://arxiv.org/abs/#1}{arXiv:#1}}

\def\cprime{$'$} \def\cprime{$'$} \def\cprime{$'$}
\providecommand{\bysame}{\leavevmode\hbox to3em{\hrulefill}\thinspace}
\providecommand{\MR}{\relax\ifhmode\unskip\space\fi MR }
\providecommand{\MRhref}[2]{%
  \href{http://www.ams.org/mathscinet-getitem?mr=#1}{#2}
}
\providecommand{\href}[2]{#2}

\end{document}